\documentclass[12pt]{amsart}
\usepackage{latexsym,amscd,amssymb,epsfig,mathrsfs} 
\usepackage[dvips]{color}
\usepackage[all]{xy}
\usepackage[T1]{CJKutf8}
\usepackage{hyperref}

\theoremstyle{plain}
  \newtheorem{theorem}{Theorem}[subsection]
  \newtheorem{proposition}[theorem]{Proposition}
  \newtheorem{lemma}[theorem]{Lemma}
  \newtheorem{corollary}[theorem]{Corollary}
  
  \newtheorem*{thma}{Theorem A}
  \newtheorem*{thmb}{Theorem B}
\theoremstyle{definition}
  \newtheorem{definition}[theorem]{Definition}
  \newtheorem{example}[theorem]{Example}

\theoremstyle{remark}
  \newtheorem{remark}[theorem]{Remark}

\numberwithin{equation}{section}

\def\umapright#1{\smash{
   \mathop{\longrightarrow}\limits^{#1}}}

\def\rmapdown#1{\Big\downarrow\rlap
   {$\vcenter{\hbox{$\scriptstyle#1$}}$}}

\def\tempbaselines
{\baselineskip22pt\lineskip3pt
   \lineskiplimit3pt}
\def\diagram#1{\null\,\vcenter{\tempbaselines
\mathsurround=0pt
    \ialign{\hfil$##$\hfil&&\quad\hfil$##$\hfil\crcr
      \mathstrut\crcr\noalign{\kern-\baselineskip}
  #1\crcr\mathstrut\crcr\noalign{\kern-\baselineskip}}}\,}

\def\pullback#1&#2&#3&#4&#5&#6&#7&#8&{
\diagram{#1&\umapright{#2}&#3\cr
\rmapdown{#4}&&\rmapdown{#5}\cr
#6&\umapright{#7}&#8\cr}}

\def\calC{{\mathcal C}}
\def\calD{{\mathcal D}}

\def\calK{{\mathcal K}}
\def\calJ{{\mathcal J}}
\def\calO{{\mathcal O}}

\def\calU{{\mathcal U}}
\def\frakM{{\mathfrak M}}
\def\frakN{{\mathfrak N}}
\def\frakF{{\mathfrak F}}
\def\frakG{{\mathfrak G}}
\def\frakO{{\mathfrak O}}
\def\frakQ{{\mathfrak Q}}

\def\k{{\underline{k}}}
\def\ZZ{\mathbb{Z}}
\def\I{\mathbb{I}}

\DeclareMathOperator{\Aut}{\mathrm{Aut}}

\DeclareMathOperator{\cH}{\check{\mathrm{H}}}

\def\colim{\mathop{\varprojlim}\nolimits}

\DeclareMathOperator{\End}{\mathrm{End}}
\DeclareMathOperator{\Ext}{\mathrm{Ext}}

\DeclareMathOperator{\GM}{\mathbb{G}_m}
\def\H{\mathop{\rm H}\nolimits}
\def\Hom{\mathop{\rm Hom}\nolimits} 
\def\Nat{\mathop{\rm Nat}\nolimits} 
 
\def\Id{\mathop{\rm Id}\nolimits}
\def\kt{\mathop{\underline{k^{\times}}}\nolimits}
\def\colim{\mathop{\varinjlim}\nolimits}
\def\lim{\mathop{\varprojlim}\nolimits}
\def\Ob{\mathop{\rm Ob}\nolimits} 
\def\Mor{\mathop{\rm Mor}\nolimits}
\def\Res{\mathop{\rm Res}\nolimits}
\DeclareMathOperator{\Pic}{Pic}
\def\PSh{\mathop{\rm PSh}\nolimits}
\def\Sh{\mathop{\rm Sh}\nolimits}

\DeclareMathOperator{\stmod}{-\underline{mod}}

\def\Stab{\mathop{\rm Stab}\nolimits}

\def\OpG{{\calO^{\circ}_p(G)}}
\def\OPG{{\mathbf{O}^{\circ}_p(G)}}

\DeclareMathOperator{\Set}{-Set}

\DeclareMathOperator{\OG}{\mathbf{O}\it{G}}

\DeclareMathOperator{\C}{\mathbf{C}}
\DeclareMathOperator{\D}{\mathbf{D}}
\DeclareMathOperator{\Mod}{-Mod}
\DeclareMathOperator{\Md}{\frakM{\rm od}-}
\DeclareMathOperator{\Ab}{\mathbf{Ab}}
\DeclareMathOperator{\rMod}{Mod-}

\begin{document}

\title[Endotrivial modules]
{On cohomological characterizations of endotrivial modules}

\author{Fei Xu}
\author{Chenyou Zheng}
\email{fxu@stu.edu.cn}
\email{18cyzheng@stu.edu.cn}

\address{Department of Mathematics\\
Shantou University\\
Shantou, Guangdong 515063, China}


\thanks{The authors \begin{CJK*}{UTF8}{}
\CJKtilde \CJKfamily{gbsn}(徐 斐、郑晨佑)
\end{CJK*} are supported by the NSFC grant No. 12171297}

\begin{abstract}
Given a general finite group $G$, there are various finite categories whose cohomology theories are of great interests. Recently Balmer and Grodal gave some new characterizations of the groups of endotrivial modules, via \v{C}ech cohomology and category cohomology, respectively, defined on certain orbit categories. These two seemingly different approaches share a common root in topos theory. We shall demonstrate the connection, which leads to a better understanding as well as new characterizations of the group of endotrivial modules.
\end{abstract}

\maketitle

\section{Introduction}

Let $G$ be a general finite group and $p$ be a prime. The category $G\Set$ of (finite left) $G$-sets and $G$-equivariant maps is a standard framework to discuss $G$-actions. Within $G\Set$, the orbit category $\calO(G)$ is a full subcategory, consisting of objects $G/H$, for all subgroups $H\subset G$. One denotes by $\calO_p(G)$ (or $\calO_p^{\circ}(G)$, respectively), the full subcategory of $\calO(G)$, consisting of orbits $G/Q$ where $Q$ is a (or a non-trivial, respectively) $p$-subgroup.

Cohomology of orbit categories plays an important role in group representations and cohomology, see \cite{Sm} and its references. We shall focus on the recent works of Balmer \cite{Ba} and Grodal \cite{Gr} on the endotrivial modules, especially on the special case of the Sylow-trivial modules. Let $k$ be a field of positive characteristic $p$ dividing the order of $G$. Suppose $k^{\times}$ is the abelian group of invertible elements in $k$. It gives rise to each of the preceding categories a constant presheaf(=contravariant functor) of abelian groups, written as $\kt$.

Consider the stable module category $kG\stmod$. Its Picard group $T_k(G):=\Pic(kG\stmod)$ is called the \textit{group of endotrivial modules}, which is an important subject in modular representation theory, see \cite{Li2, Ma}. A $kG$-module $M$ is \textit{endotrivial} if $\End_k(M)\cong k$ in $kG\stmod$, and $T_k(G)$ consists of the isomorphism classes of endotrivial modules, equipped with an operation given by the tensor product. It is strightforward to verify that $M\otimes_k-:kG\stmod\to kG\stmod$ is an equivalence if and only if $M$ is endotrivial. The classification of $T_k(G)$ has been an active direction in modular representation theory. When $G$ is a $p$-group, the classification is done \cite{Th}. While for general finite groups, it is still open.

Given a subgroup $H$ of $G$, there is a natural restriction $T_k(G)\to T_k(H)$, whose kernel is written as $T_k(G,H)$. To understand the structure of $T_k(G)$, it is key to characterize the kernel (as well as its image) when $H=P$ is a Sylow $p$-subgroup. The modules representing elements of $T_k(G,P)$ are called \textit{Sylow-trivial}. In \cite{Ba13}, Balmer identified $T_k(G,P)$ with $A_k(G,P)$, the group of weak $P$-homomorphisms. Soon after, he reinterpreted his result by a \v{C}ech cohomology group \cite{Ba} 
$$
T_k(G,P)\cong\cH^1(\mathcal{U}_{G/G},\GM).
$$ 
Recall that Balmer \cite{Ba} introduced on $G\Set$ the so-called \textit{sipp-topology}, via a specified topological basis $\calK_{\rm sipp}$, and fixed a sipp-covering $\mathcal{U}_{G/G}=\{G/P\to G/G\}\in\calK_{\rm sipp}(G/G)$. Then he took a sheaf $\GM$ on the site $(G\Set, \calK_{\rm sipp})$, determined by the constant presheaf $\kt$, to demonstrate the above isomorphism. Later on Grodal \cite{Gr} proved that (in terms of category cohomology)
$$
T_k(G,P)\cong\H^1(\calO_p^{\circ}(G),\kt).
$$

Since $\calO_p^{\circ}(G)$ is a full subcategory of $G\Set$, and $\GM$ is related to $\kt$, it is natural to ask whether these two cohomology groups match in their own right? In fact, both cohomology theories are defined on sites built upon relevant orbit categories. We shall give a general discussion about cohomology on finite sites, which eventually leads to the following general statement.

\begin{thma}[Corollary 4.4.3] Let $\frakM$ be a presheaf of abelian groups on $\OpG$. Then its right Kan extension to $\calO(G)$, along the inclusion $\iota : \OpG \to \calO(G)$, $RK_{\iota}(\frakM)$ is a sheaf with respect to the sipp topology. Suppose $\calU_{G/G}$ is the aforementioned covering of $G/G$. For any $i\ge 0$, there is a canonical isomorphism 
$$
\cH^i(\calU_{G/G}, RK_{\iota}(\frakM)) \cong \H^i(\OpG, \frakM).
$$
\end{thma}

If $\frakM=\kt$ and $i=1$, then $RK_{\iota}(\kt)\cong\GM$ and the above isomorphism provides an identification between Balmer's and Grodal's formulas on $T_k(G,P)$.

The proof of our theorem depends on a review of necessary topos theory. In fact, we shall investigate sheaves built on the following categories
$$
\calO^{\circ}_p(G)\subset\calO_p(G)\subset\calO(G)\subset G\Set,
$$
and then consider geometric morphisms between topoi (in other words, sheaf categories). The ingredients of our proof will scatter in our presentation through Sections 2, 3 and 4. Then our lengthy discussion accumulates in Section 4 to establish the above result.

Based on the same method, we are able to restate Balmer's formula for 
$$
I_k(G,P)={\rm im}(T_k(G)\to T_k(P))\cong\ker(\cH^0(\calU_{G/G},{\rm Pic})\to\cH^2(\calU_{G/G},\GM))
$$ 
with category cohomology on $\calO_p^{\circ}(G)$. By an explicit construction of the map on the right-hand-side, we may better understand its kernel.

Since both degree one \v{C}ech cohomology and category cohomology groups can be expressed in terms of topos cohomology, we use the Leray spectral sequence for a certain geometric morphism to provide another proof of $\cH^1(\calU_{G/G}, \GM) \cong \H^1(\OpG, \underline{k^{\times}})$. This approach also helps us to describe $T_k(G,P)$ itself as a Picard group. Consider the abelian category of presheaves of $k$-modules on $\calO^{\circ}_p(G)$, with the natural (objectwise) tensor product. This is the same as that of the modules on the ringed site $(\calO^{\circ}_p(G),\k)$, whose Picard group is written as $\Pic(\k)$ \cite[Definition 18.32.6]{St}. It is identified with the first cohomology of $\calO^{\circ}_p(G)$ with coefficients in $\underline{k^{\times}}$ \cite[Lemma 21.6.1]{St}.

\begin{thmb}[Theorem 6.4.1] There is an isomorphism of abelian groups
$$
T_k(G,P) \cong \Pic(\k).
$$
The group $\Pic(\underline{k})$ consists of the isomorphism classes of presheaves $\frakF$ satisfying the conditions that $\dim_k\frakF(G/Q)=1$ and that  $\frakF(G/Q\to G/Q')$ is an isomorphism of $k$-vector spaces, for every morphism $G/Q \to G/Q'$ of $\calO^{\circ}_p(G)$.
\end{thmb}

The paper is organized as follows. In Section 2, we recall basic concepts of Grothendieck sites and sheaves. Then in Section 3 we show one can work within $\calO(G)$ instead of $G\Set$, in order to rewrite Balmer's \v{C}ech cohomology group. We continue in Section 4 to establish the promised isomorphisms in Theorem A. Finally in Sections 5 and 6, we present more applications of our methods to the endotrivial modules, including Theorem B.

\section{Sheaves on sites}

For the convenience of the reader, we recall basics of Grothendieck sheaf theory here. Since we will work with finite categories and later on use Verdier's version of \v{C}ech cohomology, we will mainly follow \cite[Expos\'e i-v]{SGA4}, but some notations will be similar to the relatively modern sources \cite{Jo, MM}. To avoid repetition, we shall extract only necessary preliminary parts from \cite{WX}, with an emphasis on sites built upon the following four categories 
$$
\calO^{\circ}_p(G)\subset\calO_p(G)\subset\calO(G)\subset G\Set.
$$
Here $\calO(G)$ is the category of orbits $G/H$, with $H$ running over the set of all subgroups of $G$. It has a full subcategory $\calO_p(G)$, where $H$ is required to be a $p$-subgroup. If we discard the object $G/1$ in $\calO_p(G)$, then we obtain $\OpG$. It is well-known that the automorphism group of $G/H$ is $N_G(H)/H$.

\subsection{Topologies and sites} Let $\calC$ be a \textit{small} category. We assume the reader is familiar with the concept of a functor category. For consistency, we shall call a contravariant functor on $\calC$ a \textit{presheaf}. We denote by $\PSh(\calC)$ the category of presheaves from $\calC$ to Set, the category of sets. If $R$ is a commutative ring with identity, we write $\PSh(\calC,R)$ for the category of presheaves of $R$-modules (contravariant functors from $\calC$ to $R\Mod$). If $R=\ZZ$, this is often dubbed as $\PSh(\calC,\ZZ)=\PSh(\calC,\Ab)$.

One may put a \textit{Grothendieck topology} $\calJ_{\calC}$ on $\calC$, see \cite[II D\'efinition 1.1]{SGA4}. It relies on the concept of a \textit{sieve} $S$ on $x\in\Ob\calC$, which is simply a subfunctor of $\Hom_{\calC}(-,x)$. Moreover a sieve $S$ on $x$ is identified with a set (also written as $S$) of morphisms with codomain $x$, satisfying the condition that if ${\sf u}\in S$ and ${\sf uv}$ exists then ${\sf uv}\in S$. For instance, the \textit{maximal sieve} $\Hom_{\calC}(-,x)$ is given by the set of all morphisms with codomain $x$.

\begin{definition} Let $\calC$ be a small category. A \textit{Grothendieck topology} on $\calC$ is a function $\calJ_{\calC}$ which assigns to each object $x\in\Ob\calC$ a set of sieves $\calJ_{\calC}(x)$ on $x$, in such a way that
\begin{enumerate}
\item the maximal sieve $\Hom_{\calC}(-,x)$ is in $\calJ_{\calC}(x)$;
\item if $S\in\calJ_{\calC}(x)$, then ${\sf u}^*(S)=\{{\sf v}\bigm{|} {\sf uv}\in S\}$ lies in $\calJ_{\calC}(y)$ for any ${\sf u}:y\to x$;
\item if $S_1\in\calJ_{\calC}(x)$ and $S_2$ is any sieve on $x$ such that ${\sf u}^*(S_2)\in\calJ_{\calC}(y)$ for all ${\sf u}: y \to x$ in $S_1$, then $S_2\in\calJ_{\calC}(x)$.
\end{enumerate}
\end{definition}

Any sieve in $\calJ_{\calC}(x)$ is called a \textit{covering sieve} on $x$. A small category $\calC$ equipped with a Grothendieck topology $\calJ_{\calC}$ is called a \textit{site} $\C=(\calC,\calJ_{\calC})$ (the subscript is often omitted when there is no confusion). If the underlying category $\calC$ is finite, then we call $\C$ a \textit{finite site}.

\begin{example}
Given any category $\calC$, one can define the \textit{minimal topology} $\calJ_{\min}$ such that $\calJ_{\min}(x)=\{\Hom_{\calC}(-,x)\}, \forall x\in\Ob\calC$.

There is also a \textit{maximal topology} $\calJ_{\max}$ on $\calC$ which asks all sieves (including the empty ones) to be covering.
\end{example}

\begin{remark} It follows from the axioms that the intersection of two covering sieves on $x$ is again a covering sieve on $x$. When $\calC$ is \textit{finite}, there are always finitely many covering sieves on each $x$. Thus the intersection of these covering sieves is the \textit{minimal covering sieve} $S^{\min}_x$ on $x$ (which could be the empty sieve).

Suppose there are two sieves $S\subset R$ on $x$. If $S$ is covering, then $R$ must be covering as well. As a consequence, for a finite site $(\calC,\calJ)$, the topology is determined by the (finite) set of minimal covering sieves $\{S^{\min}_x\in\calJ(x)\bigm{|}x\in\Ob\calC\}$.
\end{remark}

\begin{remark} In practice, one often finds that the underlying category $\calC$ admits pull-backs (i.e. fibre products). Under the circumstance, it is common to introduce a \textit{topological basis} (or a pre-topologi\'e) $\calK_{\calC}$. Again, it assigns to each object $x$ a collection $\calK_{\calC}(x)$, consisting of families of morphisms with $x$ as codomian which satisfy another set of axioms \cite{MM}. We also call $\C=(\calC,\calK_{\calC})$ a site.

Each topological basis defines a unique topology $\calJ_{\calC}$ on $\calC$. Indeed, one asks a sieve $S$ on $x$ to be a covering sieve (that is, $S\in\calJ_{\calC}(x)$) if and only if there is some $R\in\calK_{\calC}(x)$ such that $R\subset S$.
\end{remark}

We present Balmer's sipp-topology (given by a topological basis) on $\calC=G\Set$ \cite{Ba}, the category of finite $G$-sets.

\begin{example} Suppose $p$ is a prime (dividing the order of $G$) and $U$ is a $G$-set. Consider a family of $G$-maps $\calU_U=\{\alpha_i : U_i \to U\}_i$. Balmer calls it a sipp-covering of $U$ if to each $u\in U$, there exists an $i_0$ and a preimage $u_{i_0}\in\alpha_{i_0}^{-1}(u)$ satisfying the condition that
$$
(p,[\Stab_G(u):\Stab_G(u_{i_0})])=1.
$$
Balmer introduced a topological basis on $G\Set$ by setting $\calK_{\rm sipp}(U)=\{\mbox{all sipp-coverings on}\ U\}, \forall U$. We should mention that ``sipp'' stands for ``stabilizer of index prime to $p$'' (or the French version, like for the old Grothendieck topologies of SGA etc. ``stabilisateur d'indice premier \'a $p$''). For instance, when $U=G/H$ and $P_H$ is a Sylow $p$-subgroup of $H$, 
$$
\calU_{G/H}=\{G/P_H \to G/H\}
$$
is a sipp-covering on $G/H$, where the morphism is given by $gP_H \mapsto gH, \forall g\in G$.

For future reference, we generate a topology $\calJ_{\rm sipp}$ on $G\Set$. It is easy to see that for a sieve $S$ on $U$ to be a covering sieve (that is $S\in\calJ_{\rm sipp}(U)$) if and only if there exists a $G$-map $\alpha: V \to U$ in $S$, matching the preceding condition. The resulting (sipp-)topology is written as $\calJ_{\rm sipp}$. For instance $\calJ_{\rm sipp}(G/H)$ is the set of sieves on $G/H$ which contains any morphism of the form $G/P_H \to G/H$.
\end{example}

There are many interesting topologies than those given in the above examples. We recall a key construction. Suppose $x\in\Ob\calC$ and $\{x_i\to x\}_i$ is a set of morphisms in $\calC$. For any $w\in\Ob\calC$, there is a natural set map
$$
\coprod_{x_i\to x}\Hom_{\calC}(w,x_i) \to \Hom_{\calC}(w,x).
$$

\begin{definition}{\cite[IV Exercise 9.1.12]{SGA4}} Let $\calC$ be a small category and $\calD$ be a strictly full subcategory. The \emph{subcategory topology} $\calJ^{\calD}$ is given by defining $\calJ^{\calD}(x)$, $\forall x\in\Ob\calC$, to consist of the sieves $S$ on $x$ which satisfy the condition that the natural set map
$$
\coprod_{y\to x \in S}\Hom_{\calC}(w,y) \to \Hom_{\calC}(w,x)
$$
is surjective for all $w \in \Ob\calD$.
\end{definition}

In Example 2.1.2, the minimal (or the maximal, respectively) topology is given by the subcategory $\calC$ (or $\emptyset$, respectively). 

\begin{lemma} Consider the category $G\Set$. 
\begin{enumerate}
\item Let $\calD$ be the full subcategory with only one object $G/1$. The corresponding subcategory topology consists of all non-empty sieves.

\item Let $\calD=\calO_p(G)$. The resulting subcategory topology identifies with Balmer's $\calJ_{\rm sipp}$. 
\end{enumerate}
\end{lemma}

\begin{proof} One just uses the definitions to verify the facts. To see (2), one may assume $x=G/H$ and compare with the characterization of $\calJ_{\rm sipp}(G/H)$ in Example 2.1.5. Suppose $S\in\calJ^{\calO_p(G)}(G/H)$. Then the natural set map
$$
\coprod_{y\to G/H \in S}\Hom_{G\Set}(w,y) \to \Hom_{G\Set}(w,G/H)
$$
is surjective for all $w=G/Q \in \Ob\calO_p(G)$ for some $p$-subgroup $Q$. Now any morphism $G/Q \to G/H$ is given by $Q \mapsto gH$ for some $g\in G$, satisfying the condition that $g^{-1}Qg\subset H$. Now we take $Q=P_H$ to be a Sylow $p$-subgroup of $H$. Then every morphism $G/P_H \to G/H$ admits a preimage in $\coprod_{y\to G/H \in S}\Hom_{G\Set}(G/P_H,y)$. It means that all $G/P_H \to G/H$ lie in $S$. Therefore $\calJ^{\calO_p(G)}(G/H)=\calJ_{\rm sipp}(G/H)$, for all $G/H$ (and thus every $x\in\Ob G\Set$).

The proof of (1) is similar. This topology on $G\Set$ is called the atomic topology $\calJ_{\rm at}$, studied by M. Artin \cite{MM}. 
\end{proof}

A sieve like $S^{\min}_{G/H}$ is called a {\it principal sieve}, as it consists of all morphisms $ff'$, where $f'$ runs over the set of all morphisms that is composable with $f$ from the right. If $H$ is a $p$-group, $S^{\min}_{G/H}$ equals the maximal sieve $h_{G/H}$, generated by any isomorphism (e.g. the identity map) $G/H \to G/H$. If $p \nmid |H|$, $S^{\min}_{G/H}$ is generated by any $G$-map $G/1\to G/H$. Balmer \cite{Ba} considered the sipp-covering of the terminal object $G/G$
$$
\calU_{G/G}=\{G/P\to G/G\}\in\calK_{\rm sipp}(G/G).
$$
It generates the minimal covering sieve $S^{\min}_{G/G}\in\calJ_{\rm sipp}(G/G)$.

We shall see that Balmer's constructions on $G\Set$ can be reduced to a site on $\calO(G)$, without losing any information.

\subsection{Sheaves} Since we will speak about \v{C}ech cohomology, as well as cohomology of topoi, we recall the definition of a sheaf on a given site $\C=(\calC,\calJ)$. Roughly speaking, a sheaf of sets on $\calC$ is a presheaf of sets on $\calC$, satisfying certain glueing properties mandated by $\calJ_{\calC}$. A concise definition of a sheaf is the following.

\begin{definition}
A presheaf $\frakF\in\PSh(\calC)$ is a ($\calJ$-)sheaf if, for every $x\in\Ob\calC$ and every $S\in\calJ(x)$, the inclusion $S\hookrightarrow\Hom_{\calC}(-,x)$ induces an isomorphism
$$
\Nat(\Hom_{\calC}(-,x),\frakF){\buildrel{\cong}\over{\to}}\Nat(S,\frakF).
$$
The category of sheaves of sets on $\C=(\calC,\calJ)$ is the full subcategory of $\PSh(\calC)$, consisting all all ($\calJ$-)sheaves, denoted by $\Sh(\C)=\Sh(\calC,\calJ)$.

A \textit{(Grothendieck) topos} $\mathscr{E}$ is a category that is equivalent to some $\Sh(\C)$, a sheaf category of sets on a site $\C$. 
\end{definition}

Note that by Yoneda Lemma, $\frakF(x)\cong\Nat(\Hom_{\calC}(-,x),\frakF)$. One may define sheaves on a site $(\calC,\calK)$ as well. If $\calJ$ is the topology generated by a topological basis $\calK$, then there is a canonical identification of sheaves. In short, one has $\Sh(\calC,\calK)=\Sh(\calC,\calJ)$. 

\begin{example}
Given any category $\calC$, one can see $\Sh(\calC,\calJ_{\min})=\PSh(\calC)$ and $\Sh(\calC,\calJ_{\max})$ is the trivial category.
\end{example}

\begin{remark} To find applications in representation theory and to introduce cohomology theory, one considers sheaves with algebraic structures, such as sheaves of $R$-modules for a given ring $R$ and sheaves of rings. By definition, a presheaf of $R$-modules is a sheaf if and only if, considered as a presheaf of sets, it is a sheaf. Sheaves of $R$-modules on a site $\C$ form a category $\Sh(\C,R)$, see \cite[7.44]{St}, which is abelian. This category has enough injectives, but usually not enough projectives.

The objects of $\Sh(\C,R)$ are exactly the $R$\textit{-module objects} in $\Sh(\C)$. Especially, the objects of $\Sh(\C,\mathbb{Z})=\Sh(\C,\Ab)$ are the abelian group objects in $\Sh(\C)$, that is, $Ab(\Sh(\C))=\Sh(\C,\Ab)$.

As a strategy, one can first prove results about topoi (categories of sheaves of sets). Then there are ways to deduce corresponding statements about sheaves with algebraic structures \cite[7.43 \& 18.13]{St}.
\end{remark}

\subsection{Picard groups} Let $(\calC,\otimes,\mathbf{1})$ be a symmetric monoidal category. The following construction is known to the experts.

\begin{definition} Given a symmetric monoidal category 
$(\calC,\otimes,\mathbf{1})$, an object $x\in\Ob\calC$ is said to be invertible if there exists another object $y$ such that $x\otimes y\cong\mathbf{1}$. The Picard group $\Pic(\calC)=\Pic(\calC,\otimes,\mathbf{1})$ is the group of isomorphism classes of invertible objects, with respect to $\otimes$.
\end{definition}

We shall be interested in the following two cases.

\begin{example} 
\begin{enumerate} 
\item Let $G$ be a finite group and $k$ a field. The stable module category $kG\stmod$ forms a symmetric monoidal category $(kG\stmod,\otimes_k,k)$. Its Picard group $\Pic(kG\stmod)$ is exactly $T_k(G)$.

\item Let $\C$ be a site with a structure sheaf $\frakO$ (a sheaf of commutative rings). Then the category of modules $\Md\frakO$ on the ringed site $(\C,\frakO)$ leads to a symmetric monoidal category $(\Md\frakO,\otimes,\frakO)$. Its Picard group $\Pic(\Md\frakO)$ is often abbreviated as $\Pic(\frakO)$. We mainly focus on the constant structure sheaves $\underline{k}$, in which case $\Md\underline{k}$ is the same as $\Sh(\C,k)$, the category of sheaves of $k$-modules.
\end{enumerate} 
\end{example}

\section{Morphisms of topoi and reduction on sheaf categories} 

Recall that we have the following four categories
$$
\calO^{\circ}_p(G)\subset\calO_p(G)\subset\calO(G)\subset G\Set.
$$
Balmer worked on $G\Set$, while Grodal dealt with $\calO^{\circ}_p(G)$. To establish connections between their formulas, we need two steps of reductions: firstly, we reduce from sipp-site on $G\Set$ to one on $\calO(G)$, and secondly we reduce from $\calO(G)$ to $\calO^{\circ}_p(G)$. Both steps depends on morphisms of sites and topoi. However, these two stages use different techniques, as we will see in Sections 3 and 4. (It is possible to use $\calO_p(G)$ as the intermediate category for reductions, but it needs an extra level of technicality from topos theory.)

For future reference, we recall some fundamental concepts and constructions, for comparing sheaves, see for instance \cite{SGA4, Jo, MM}. Given a (covariant) functor $\alpha:\calD\to\calC$, there exists a \textit{restriction} along $\alpha$, $Res_{\alpha} : \PSh(\calC)\to\PSh(\calD)$ which is given by the precomposition with $\alpha$. This functor admits two adjoint functors, the left and right \textit{Kan extensions} along $\alpha$, written as $LK_{\alpha}, RK_{\alpha} : \PSh(\calD)\to\PSh(\calC)$. The \textit{comma categories} (including undercategories and overcategories) are used to define the Kan extensions. Let $x\in\Ob\calC$. Then the \textit{category over} $x$, or just \textit{overcategory}, $\alpha/x$ has objects $\{(d,{\sf t})\bigm{|} {\sf t}: \alpha(d) \to x, d\in\Ob\calD\}$. A morphism $(d, {\sf t})\to(d', {\sf t}')$ is given by ${\sf u}:d\to d'$ such that ${\sf t}'\alpha({\sf u})={\sf t}$.

\subsection{Morphisms of topoi}

To compare sheaves on two sites, we use the concepts of continuous and cocontinuous functors \cite[Expos\'e III]{SGA4}. Recall that given a set of morphisms $\{{\sf u}_i\}_{i\in I}$ with common codoamin $x\in\Ob\calC$, the sieve $S$ generated by $\{{\sf u}_i\}_{i\in I}$ is the smallest sieve that contains these ${\sf u}_i$'s, written as $S=({\sf u}_i | i\in I)$.

\begin{definition} Let $\C=(\calC,\calJ_{\calC})$ and $\D=(\calD,\calJ_{\calD})$ be two sites. A functor $\alpha: \calD\to\calC$ is called \textit{continuous} if, for each $S_d\in\calJ_{\calD}(d)$, $\alpha(S_d)$ generates a covering sieve on $\alpha(d)$ and, for each $x\in\Ob\calC$, $(x/\alpha)^{op}\simeq\alpha^{op}/x^{op}$ is filtered. 

A functor $\beta : \calD\to\calC$ is called \textit{cocontinuous} if for any covering sieve $S_{\beta(d)}\in\calJ_{\calC}(\beta(d))$ on $\beta(d)\in\Ob\calC$, there exists a covering sieve $S_d\in\calJ_{\calD}(d)$ on $d\in\Ob\calD$ with $\beta(S_d)\subset S_{\beta(d)}$. (This is equivalent to saying that $\{{\sf u}\in\Hom_{\calD}(-,d)\bigm{|}\beta({\sf u})\in S_{\beta(d)}\}$ is a covering sieve on $d$.)
\end{definition}

A functor $\alpha$ between two sites being continuous is equivalent to asking $Res_{\alpha}$ to preserve sheaves. While it being cocontinuous is equivalent to saying that $RK_{\alpha}$ preserves sheaves. We take the above forms because they are easy to check.

\begin{lemma} Consider the site $(G\Set,\calJ_{\rm sipp})$, and its subsites $\OG:=(\calO(G),\calJ_{\rm sipp})$ and $\OPG:=(\OpG,\calJ_{\min})$. The inclusion functors $\OpG \to \calO(G)$ and $\calO(G) \to G\Set$ are both continuous and cocontinuous.
\end{lemma}

\begin{proof} Verify with the definitions.
\end{proof}

The inclusion functors are encoded into geometric morphisms of topoi that we will come back to later on.

\begin{definition} A \textit{(geometric) morphism of topoi} 
$$
\Psi=(\Psi^{-1},\Psi_*) : \Sh(\D)\to\Sh(\C)
$$ 
consists of a pair of functors $\Psi_* : \Sh(\D)\to\Sh(\C)$ and $\Psi^{-1}:\Sh(\C)\to\Sh(\D)$ such that $\Psi^{-1}$ is left exact and is left adjoint to $\Psi_*$.
\end{definition}

Either a continuous or a co-continuous functor induces a morphism of topoi. We record the following (\cite[IV 4.7]{SGA4}, also see \cite[7.21.5]{St} and the comments before \cite[C.2.3.23]{Jo}).

\begin{proposition} Let $\C$ and $\D$ be two sites. If there is a continuous and cocontinuous functor $\phi:\calD\to\calC$, we have the associated morphism of topoi $\Phi=(\Phi^{-1},\Phi_*) : \Sh(\D)\to\Sh(\C)$, 
\begin{enumerate}
\item $\Phi^{-1}=Res_{\phi}$,
\item $\Phi_*=RK_{\phi}$,
\item $\Phi_!=LK_{\phi}^{\sharp}$ is left adjoint to $\Phi^{-1}$.
\end{enumerate}
\end{proposition}

Here for a sheaf $\frakG$ on $\D$, $LK_{\phi}^{\sharp}\frakG:=(LK_{\phi}\frakG)^{\sharp}$, where $(-)^\sharp : \PSh(\calC) \to \Sh(\C)$ is \textit{sheafification}. It will be discussed shortly in Section 4.2.

Many results on topoi may be readily applied to sheaves with algebraic structures. This is vital for investigating representations and cohomology. From \cite[7.43]{St} the aforementioned constructions and properties we know about sheaves of sets pass to sheaves of abelian groups, $R$-modules etc., in the sense that the morphism $\Psi$ of topoi gives rise to a pair of adjoint functors $(\Psi^{-1},\Psi_*)$ on categories of sheaves with suitable algebraic structures.

\subsection{The Comparison Lemma and reduction to $\calO(G)$}

The inclusion $\calO(G) \to G\Set$ has rather nice properties, which allow us to invoke the Comparison Lemma. Then we can replace Balmer's $(G\Set,\calJ_{\rm sipp})$ by an appropriate site on the full subcategory $\calO(G)$. We follow the presentation of the Comparison Lemma in \cite[C.2.?]{Jo}.

\begin{definition} Let $\C=(\calC,\calJ)$ be a site and $\calD\subset\calC$ be a full subcategory. Define a topology $\calJ|_{\calD}$ on $\calD$ by $\calJ|_{\calD}(y)=\{S\cap\Mor\calD\bigm{|} S\in\calJ(y)\}, \forall y\in\Ob\calD$. This is the \textit{induced topology} on $\calD$.

We call $\D=(\calD,\calJ|_{\calD})$ a \textit{subsite} of $\C=(\calC,\calJ)$.
\end{definition}

\begin{definition} Let $\C=(\calC,\calJ_{\calC})$ be a site and $\calD\subset\calC$ be a full subcategory. If, for every $x\in\Ob\calC$, there exists a covering sieve $S\in\calJ_{\calC}(x)$ such that ${\rm dom}(\alpha)\in\Ob\calD, \forall \alpha\in S$, then we call $\D$ a \textit{dense subsite} of $\C$, where $\calD$ is given the induced topology.
\end{definition}

\begin{proposition}[Comparison Lemma] If $\D$ is a dense subsite of $\C$, then the inclusion functor $\iota: \calD \to \calC$ induces an equivalence  $\Sh(\D)\simeq\Sh(\C)$.
\end{proposition}

The equivalences are given by the restriction and the right Kan extension, along the inclusion $\calD \hookrightarrow \calC$, as in Proposition 3.1.4.

\begin{example} Let $(G\mbox{-}\rm{Set}, \calJ)$ be a site. Then the induced subsite $(\calO(G),\calJ|_{\calO(G)})$ is always dense. Therefore 
$$
\Sh(G\mbox{-}\rm{Set}, \calJ)\simeq\Sh(\calO(G),\calJ|_{\calO(G)}).
$$
\end{example}

\begin{example} It is known that each topology $\calJ$ on a finite category $\calC$ is uniquely given by a strictly full idempotent-complete subcategory $\calD$ and is identified with $\calJ^{\calD}$. Moreover, $(\calD,\calJ^\calD|_\calD)$ is a dense subsite of $(\calC,\calJ^\calD)$, with $\calJ^\calD|_\calD=\calJ_{\min}$ the minimal topology on $\calD$. It follows from the Comparison Lemma that
$$
\Sh(\calC,\calJ^\calD)\simeq\PSh(\calD).
$$
\end{example}

Recall that $\calO(G)$ is the category of orbits $G/H$, where $H$ runs over the set of all subgroups of $G$. It has a full subcategory $\calO_p(G)$, where $H$ is required to be a $p$-subgroup. If we discard the object $G/1$ in $\calO_p(G)$, then we obtain $\OpG$. The functors in Lemma 3.1.2 leads to the following equivalences.

\begin{corollary} There are equivalences
$$
\Sh(G\Set,\calJ_{\rm sipp})\simeq\Sh(\calO(G),\calJ_{\rm sipp})\simeq\Sh(\calO_p(G),\calJ_{\min})=\PSh(\calO_p(G)).
$$
\end{corollary}

\begin{proof} The first equivalence follows from Example 3.2.4 (For brevity, we simplify $\calJ_{\rm sipp}|_{\calO(G)}$ to $\calJ_{\rm sipp}$). 

If we continue to examine the subcategory topology on $\calO(G)$, given by $\calO_p(G)$, then it is exactly $\calJ_{\rm sipp}$, and the second equivalence comes from Example 3.2.5.
\end{proof}

\begin{remark} 
\begin{enumerate}
\item The above corollary permits us to reduce from $G\Set$ to $\calO_p(G)$, better than we've claimed. However, since $\calO(G)$ has a terminal object $G/G$, and $\calO_p(G)$ does not have one in general, it is more convenient to use $\OG$ as an intermediate site to carry on our reduction in the next stage.

\item Further restriction down to the full subcategory $\calO_p^{\circ}(G)$ with induced topology does \textit{not} produce a dense subsite of $\OG=(\calO(G),\calJ_{\rm sipp})$. We cannot obtain a reduction on the level of sheaf categories. So we turn to simplify cohomology instead in the next stage.
\end{enumerate}
\end{remark}

To enter the second stage of reduction, we need to recall some other constructions on various cohomology theories on sites and topoi, in upcoming section. We end this section with a useful remark.

\begin{remark} Since $G/1$ belongs to $\calO_p(G)$, the atomic topology ($\calJ_{\rm at}=\calJ^{G/1}$) is larger than the sipp-topology ($\calJ_{\rm sipp}=\calJ^{\calO_p(G)}$, see Lemma 2.1.7), in the sense that $\calJ_{\rm sipp}\subset\calJ_{\rm at}$. Equivalently, it means that 
$$
G\Set\simeq\Sh(\calO(G),\calJ_{\rm at})\subset\Sh(\calO(G),\calJ_{\rm sipp})\simeq\PSh(\calO_p(G)).
$$
It is known \cite{MM} that the left equivalence is given by $M \mapsto \frakF_M$, the \textit{fixed-point sheaf}, which is defined on objects by $\frakF_M(G/H)=M^H$. The above inclusion implies that $\frakF_M\in\Sh(\calO(G),\calJ_{\rm sipp})$. It corresponds to its restriction, along $\calO_p(G) \subset \calO(G)$, in $\PSh(\calO_p(G))$.
\end{remark} 

\section{Reduction on \v{C}ech cohomology}

We will recall \v{C}ech cohomology here \cite[Expos\'e V, 2.4.3]{SGA4}. Since it is derived from the {\it half-sheafification} functor $(-)^{\dagger} : \PSh(\calC) \to \PSh(\calC)$, we shall present both of them. Then on $\C=\OG$ we compute relevant \v{C}ech cohomology groups.

\subsection{\v{C}ech cohomology}

\begin{definition} Let $\C=(\calC,\calJ)$ be a site and $x\in\Ob\calC$. Given $S\in\calJ(x)$, the right derived functors of $\Nat(S,-): \PSh(\calC,\Ab)\to\Ab$ are written as $\cH^i(S,-), i\ge 0$. If $\frakM\in\PSh(\calC,\Ab)$, then $\cH^i(S,\frakM), i\ge 0,$ are the {\it \v{C}ech cohomology groups} of $S$ with coefficients in $\frakM$.
\end{definition}

If $\calC$ admits pull-backs and a topological basis $\calK$ that generates $\calJ$, specifically with $\calU\in\calK(x)$ giving rise to $S\in\calJ(x)$, then one normally uses a simplicial complex to define $\cH^i(\calU,\frakM)$, which is canonically isomorphic to $\cH^i(S,\frakM), \forall i\ge 0$.

\begin{definition} Let $\C=(\calC,\calJ)$ be a site, $x\in\Ob\calC$ and $\frakM\in\PSh(\calC,\Ab)$. Then we call the following
$$
\cH^i(x,\frakM):=\colim_{\calJ(x)}\cH^i(-,\frakM),
$$
 the {\it \v{C}ech cohomology groups} of $x$ with coefficients in $\frakM$, $\forall i \ge 0$. If $\calC$ admits a terminal object $x_0$, then we call $\cH^i(\C,\frakM):=\cH^i(x_0,\frakM)$ the \v{C}ech cohomology groups of $\C$ with coefficients in $\frakM$.
\end{definition}

\begin{theorem}
Let $\C=(\calC,\calJ)$ be a site such that $\calC$ is a finite category. For every $x\in\Ob\calC$, there exists the smallest covering sieve $S^{\min}_x\in\calJ(x)$. The \v{C}ech cohomology groups can be computed by
$$
\cH^i(x,\frakM):=\colim_{\calJ(x)}\cH^i(-,\frakM)\cong\cH^i(S^{\min}_x,\frakM).
$$
\end{theorem}

\begin{proof} By Remark 2.1.3, $\calJ(x)$ has an initial object $S^{\min}_x$. Our identification of the cohomology groups follows from this, and the definition of the direct limit.
\end{proof}

With Example 2.1.5 and Lemma 2.1.7, we now specialize to the site $\OG=(\calO(G),\calJ_{\rm sipp})$.

\begin{corollary} For any $\frakM\in\PSh(\calO(G),\Ab)$ and any $i\ge 0$, 
$$
\cH^i(\calU_{G/H},\frakM)\cong\cH^i(S^{\min}_{G/H},\frakM).
$$
Particularly
$$
\cH^i(\calU_{G/G},\frakM)\cong\cH^i(S^{\min}_{G/G},\frakM)\cong\cH^i(G/G,\frakM)=\cH^i(\OG,\frakM).
$$
\end{corollary}

We shall compute $\cH^i(\calU_{G/G},\frakM)$, and then identify it with Grodal's formula.

\subsection{Sheafification} The abelian presheaf $\cH^0(-,\frakM)$ is also written as $\frakM^{\dagger}$, sometimes called the {\it half-sheafification} of $\frakM$, because $\frakM^{\sharp}=(\frakM^{\dagger})^{\dagger}$ is always an abelian sheaf on $\C$. The {\it sheafification} functor 
$$
(-)^{\sharp} : \PSh(\calC,\Ab)\to\Sh(\C,\Ab)
$$ 
is exact, left adjoint to the forgetful functor $\Sh(\C,\Ab)\to\PSh(\calC,\Ab)$, which is left exact. In fact, the formula 
$$
\frakM^{\dagger}(x)=\colim_{S\in\calJ(x)}\Nat(-,\frakM)
$$ 
works for presheaves of sets! There are two functors $(-)^\dagger : \PSh(\calC) \to \PSh(\calC)$ and $(-)^{\sharp} : \PSh(\calC) \to \Sh(\C)$.

\begin{proposition} The constant presheaf $\kt\in\PSh(\calO(G))$ is a sheaf in $\Sh(\OG)$.
\end{proposition}

\begin{proof} For each $G/H\in\Ob\calO(G)$, we have
$$
(\kt)^{\dagger}(G/H)=\colim_{\calJ_{\rm sipp}(G/H)}\Nat(-,\kt)\cong\Nat(S^{\min}_{G/H},\kt)\cong k^{\times}.
$$
One can continue to verify that $(\kt)^{\dagger}\cong\kt$. Thus $\kt\cong(\kt)^{\sharp}$ is already a sheaf.
\end{proof}

\subsection{The sheaf $\GM$} Now we turn to the sheaf $\GM$ used by Balmer, as the coefficients for \v{C}ech cohomology. We indicate that it is closely related to $\underline{k^{\times}}$ used by Grodal, as the coefficients for category cohomology.

For the reader's convenience, we apply Proposition 3.1.4 to the functor $\iota:\OpG\to\calO(G)$ and state the consequences explicitly.

\begin{corollary} Consider the two sites $\OG$ and $\OPG$. The inclusion functor $\iota:\OpG \to \calO(G)$ gives rise to a morphism of topoi $\I=(\I^{-1},\I_*): \Sh(\OPG)=\PSh(\OpG) \to \Sh(\OG)$, with
\begin{enumerate}
\item $\I^{-1}=Res_{\iota}$,
\item $\I_*=RK_{\iota}$,
\item $\I_!=LK_{\iota}^{\sharp}$ is left adjoint to $\I^{-1}$.
\end{enumerate}
\end{corollary}

It means the right Kan extenson of every presheaf on $\OpG$ is a sheaf on $\OG$. Now we can produce the sheaf $\GM$ used by Balmer \cite{Ba}.

\begin{proposition} From $\kt\in\PSh(\OpG,\Ab)$, we obtain a sheaf $\GM=\I_*(\kt)\in\Sh(\OG,\Ab)$. It is characterized by $\GM(G/H)=k^{\times}$ if $p\bigm{|} |H|$, and zero otherwise. (Here zero is understood as the terminal object of $\Ab$.)
\end{proposition}

\begin{proof} By Corollary 4.3.1, $\GM$ is a sheaf. Now by the definition of the right Kan extension, we need to understand the structures of $\iota/(G/H)$ for all $G/H$. One can verify that $\iota/(G/H)$ is connected if $p\bigm{|}|H|$, and empty otherwise. Thus $\GM(G/H)=\lim_{\iota/(G/H)}\kt$ gives the above characterization. (Here a category is said to be connected if one can find a zig-zag of some morphisms between any two objects.)
\end{proof}

We indicate that $\GM$ is not the sheafification of $\kt\in\PSh(\calO(G),\Ab)$, a minor error stated in \cite{Ba}, see Proposition 4.2.1.

After preparation on modules and cohomology on (finite) sites, we are at the position to rewrite Balmer's formulas, and compare one of them with Grodal's work. On the way, we obtain some other interesting results on endotrivial modules.

We refer to \cite{We} for basics of category cohomology.

\subsection{The kernel $T_k(G,P)$}

\begin{lemma} The minimal covering sieve $S^{\min}_{G/G}\in\PSh(\calO(G))$ restricts to a terminal object of $\PSh(\calO_p(G))$, and of $\PSh(\OpG)$.
\end{lemma}

\begin{proof} This based on our description of $S^{\min}_{G/G}(G/H)$. Since $G/G$ is the terminal object of $\calO(G)$, there is exactly one $G$-map $G/H \to G/G$. When $H$ is a $p$-group, $S^{\min}_{G/G}(G/H)$ is the set of this unique $G$-map. Otherwise $S^{\min}_{G/G}(G/H)=\emptyset$.
\end{proof}

Let $\mathbb{Z}S$ be the free abelian group object generated by a sieve $S$. Then $\mathbb{Z}S^{\min}_{G/G}\in\PSh(\calO(G),\Ab)$ is characterized by $\mathbb{Z}S^{\min}_{G/G}(G/H)\cong\mathbb{Z}$ if $H$ is a $p$-group, and $\mathbb{Z}S^{\min}_{G/G}(G/H)=0$ otherwise.

\begin{proposition} Let $\calU$ be a sipp-covering of $G/H$ and $S$ be the sieve generated by $\calU$. For any $\frakM\in\PSh(\calO(G),\Ab)$ and any $i\ge 0$, 
$$
\cH^i(\calU,\frakM)\cong\cH^i(S,\frakM)\cong\Ext^i_{\mathbb{Z}\calO(G)}(\mathbb{Z}S,\frakM).
$$
\end{proposition}

\begin{proof} The second isomorphism comes from the adjunction
$$
\Nat_{\PSh(\calO(G))}(S,\frakM)\cong\Nat_{\PSh(\calO(G),\Ab)}(\mathbb{Z}S,\frakM).
$$
We note that $\PSh(\calO(G),\Ab)\simeq\rMod\mathbb{Z}\calO(G)$ \cite{We}.
\end{proof}

When $\calU=\calU_{G/G}$ and $\frakM=\GM$, we can further reduce the calculations to $\calO^{\circ}_p(G)$.

\begin{corollary} We have 
$$
\cH^i(\calU_{G/G},\GM)\cong\cH^i(S^{\min}_{G/G},\GM)\cong\Ext^i_{\mathbb{Z}\OpG}(\underline{\mathbb{Z}},\kt).
$$
\end{corollary}

\begin{proof} It follows from \cite[Proposition 4.2.1]{Xu}. Or we just directly take a projective resolution of $\mathbb{Z}S^{\min}_{G/G}$, formed of representable functors, supported on $G/H$ where $H$ is a $p$-group. Using the fact that $\GM$ is supported on the non-trivial $p$-groups and their overgroups, we obtain, for all $i\ge 0$,
$$
\Ext^i_{\mathbb{Z}\calO(G)}(\mathbb{Z}S^{\min}_{G/G},\GM)\cong\Ext^i_{\mathbb{Z}\OpG}(\mathbb{Z}S^{\min}_{G/G}\bigm{|}_{\OpG},\GM\bigm{|}_{\OpG})
$$
The latter is exactly $\Ext^i_{\mathbb{Z}\OpG}(\underline{\mathbb{Z}},\kt)$.
\end{proof}

The canonical isomorphism $\Ext^i_{\mathbb{Z}\OpG}(\underline{\mathbb{Z}},\kt)\cong\H^i(\OpG,\kt)$ implies that Balmer's and Grodal's formulas for $T_k(G,P)$ (when $i=1$) are naturally identified.

\section{Further applications to endotrivial modules}

Our methods give rise to more information on endotrivial modules. Here we present some new characterizations of $T_k(G,P)$ and $I_k(G,P)$.

\subsection{On $T_k(G,P)$}

Consider the short exact sequence of $\mathbb{Z}\calO(G)$-modules (that is, abelian presheaves on $\calO(G)$)
$$
0 \to \mathbb{Z}S^{\min}_{G/G} \to \underline{\mathbb{Z}} \to \frakQ \to 0,
$$
where $\frakQ$ is the quotient module. It gives rise to a long exact sequence
$$
0 \to \Hom_{\mathbb{Z}\calO(G)}(\frakQ,\Pic) \to \Hom_{\mathbb{Z}\calO(G)}(\underline{\mathbb{Z}},\Pic) \to \Hom_{\mathbb{Z}\calO(G)}(\mathbb{Z}S^{\min}_{G/G},\Pic) 
$$
$$
\to \Ext^1_{\mathbb{Z}\calO(G)}(\frakQ,\Pic) \to \Ext^1_{\mathbb{Z}\calO(G)}(\underline{\mathbb{Z}},\Pic) \to \Ext^1_{\mathbb{Z}\calO(G)}(\mathbb{Z}S^{\min}_{G/G},\Pic)
$$
Since $\calO(G)$ has a terminal object $G/G$, $T_k(G)\cong\Hom_{\mathbb{Z}\calO(G)}(\mathbb{Z},\Pic)$, the restriction $T_k(G) \to T_k(P)$ can be realized as 
$$
\Hom_{\mathbb{Z}\calO(G)}(\mathbb{Z},\Pic) \to \Hom_{\mathbb{Z}\calO(G)}(\mathbb{Z}S^{\min}_{G/G},\Pic) =\Hom_{\mathbb{Z}\OpG}(\mathbb{Z},\Pic),
$$
where $\Hom_{\mathbb{Z}\OpG}(\mathbb{Z},\Pic)\subset T_k(P)$.

\begin{proposition} With the above notations, we have 
$$
T_k(G,P)\cong\Hom_{\mathbb{Z}\calO(G)}(\frakQ,\Pic)
$$ 
and 
$$
I_k(G,P)\cong\ker\{\Hom_{\mathbb{Z}\OpG}(\mathbb{Z},\Pic)\to\Ext^1_{\mathbb{Z}\calO(G)}(\frakQ,\Pic)\}.
$$
\end{proposition}

Note that based on the same techniques in proving Corollary 5.1.3 $\Ext^1_{\mathbb{Z}\calO(G)}(\frakQ,\Pic)\cong\Ext^1_{\mathbb{Z}\overline{\OpG}}(\underline{\mathbb{Z}},\Pic)$, where $\overline{\OpG}$ consists of objects $G/H$ such that $([G:H],p)=1$.

\subsection{The image $I_k(G,P)$}

In \cite{Ba}, Balmer also provided a formula 
$$
I_k(G,P)=\ker(\cH^0(\calU_{G/G}, \Pic){\buildrel{z}\over{\to}}\cH^2(\calU_{G/G},\GM)).
$$
In the same fashion, we have 
$$
\cH^0(\calU_{G/G},{\rm Pic})\cong\Hom_{\mathbb{Z}\OpG}(\underline{\mathbb{Z}},{\rm Pic})\cong\lim_{G/Q\in\OpG}T_k(Q),
$$ 
isomorphic to $\Hom_{\mathbb{Z}[N_G(P)/P]}(\mathbb{Z},{\rm Pic}(G/P))\cong T_k(P)^{N_G(P)/P}$, and 
$$
\cH^2(\calU_{G/G},\GM)\cong\Ext^2_{\mathbb{Z}\OpG}(\underline{\mathbb{Z}},\kt).
$$
We want to construct Balmer's map in terms of category cohomology
$$
z : \Hom_{\mathbb{Z}\OpG}(\underline{\mathbb{Z}},{\rm Pic}) \to \H^2(\OpG,\kt)=\Ext^2_{\mathbb{Z}\OpG}(\underline{\mathbb{Z}},\kt).
$$
Before we doing so, we must warn the reader that the abelian group $k^{\times}$ is multiplicative. Thus one has to adapt the definition in \cite{We} to our situation.

Any element $\theta\in\Hom_{\mathbb{Z}\OpG}(\underline{\mathbb{Z}},{\rm Pic})$ specifies a compatible set of elements of $\{\theta_{G/Q}(1)\in\Pic(kQ\stmod)\}_Q$, where $Q$ runs over the set of all non-trivial $p$-subgroups. Let $\theta_{G/P}(1)=[N]\in\Pic(kP\stmod)$. A representative $N\in kP\stmod$ comes with a set of isomorphisms 
$$
\{\phi_f:N\downarrow_Q \to r_f(N\downarrow_{Q'})\}_{f:G/Q\to G/Q'},
$$
where $r_f$ is the restriction along $f$. The above is a set of descent data. By Balmer's theorem \cite{Ba}, there exists some $M\in kG\stmod$ such that $M\downarrow_Q\cong N\downarrow_Q$, along with compatible morphisms. From the detection theorem \cite[Theorem 2.1]{Ma}, $M$ must be endotrivial. Without loss of generality, we assume $N=M\downarrow_P$. Note that the automorphism group of $M\downarrow_Q\in kQ\stmod$ is $k^{\times}$. Each $\phi_f$ belongs to $\Hom_{kQ\stmod}(M\downarrow_Q, r_f(M\downarrow_{Q'}))\cong k^{\times}$. 

We define a presheaf $\frakF_{\theta}: \OpG \to \Ab$ such that $\frakF_{\theta}(G/Q)=\Aut_{kQ\stmod}(M\downarrow_Q)$ and such that
$$
\frakF_{\theta}(G/Q {\buildrel{f}\over{\to}} G/Q'): \Aut_{kQ'\stmod}(M\downarrow_{Q'})\to\Aut_{kQ\stmod}(M\downarrow_Q)
$$ 
is given by $t\mapsto \phi_f^{-1}r_f(t)\phi_f$. Since $\Aut_{kQ'\stmod}(M\downarrow_{Q'})\cong k^{\times}$, the map $\frakF(G/Q {\buildrel{f}\over{\to}} G/Q')$ corresponds to the identity $1_{k^{\times}}$ and $\frakF_{\theta}\cong\underline{k^{\times}}$.

In this way, $\theta$ leads to an element in 
$$
C^2(\OpG;\frakF_{\theta})=\prod_{{\scriptscriptstyle G/Q {\buildrel{f_1}\over{\to}} G/Q' {\buildrel{f_2}\over{\to}} G/Q''}}\frakF_{\theta}(G/Q).
$$
The above constructions are independent from the choice of $M$, and thus altogether provides a map 
$$
\tilde{z} : \Hom_{\mathbb{Z}\OpG}(\underline{\mathbb{Z}},{\rm Pic}) \to \prod_{{\scriptscriptstyle G/Q {\buildrel{f_1}\over{\to}} G/Q' {\buildrel{f_2}\over{\to}} G/Q''}}\frakF_{\theta}(G/Q)
$$
by $\tilde{z}(\theta)=\{\phi_{1_{G/Q}}\in\frakF_{\theta}(G/Q)\}_{f:G/Q{\buildrel{f_2}\over{\to}} G/Q'{\buildrel{f_1}\over{\to}} G/Q''}$.

Taking into account the differential
$$
\partial_3 : \prod_{{\scriptscriptstyle G/Q \to G/Q' \to G/Q''}}\frakF_{\theta}(G/Q) \to \prod_{{\scriptscriptstyle G/R \to G/R' \to G/R''\to G/R'''}}\frakF_{\theta}(G/R),
$$
one deduces from the fact that $\theta$ is a natural transformation
$$
\begin{array}{ll}
& \partial_3(\tilde{z}(\theta))_{{\scriptscriptstyle G/R {\buildrel{f_3}\over{\to}} G/R' {\buildrel{f_2}\over{\to}} G/R''{\buildrel{f_1}\over{\to}} G/R'''}} \\
 = & \frakF(G/R{\buildrel{f_3}\over{\to}} G/R')(\phi_{1_{G/R'}})\phi_{1_{G/R}}^{-1}\phi_{1_{G/R}}\phi_{1_{G/R}}^{-1}\\
= & [\phi_{f_3}^{-1}r_{f_3}(\phi_{1_{G/R'}})\phi_{f_3}]\phi^{-1}_{G/R}\\
= & \phi_{1_{G/R}}\phi_{1_{G/R}}^{-1}\\
= & \Id_{M\downarrow_R}.
\end{array}
$$
The third equality depends on the fact that $\phi_f$'s come from a set of descent data. 

In summary, for all $\theta\in\Hom_{\mathbb{Z}\OpG}(\underline{\mathbb{Z}},{\rm Pic})$, $\tilde{z}(\theta)\in\ker(\partial_3)$ and thus $\tilde{z}$ induces Balmer's homomorphism
$$
z : \Hom_{\mathbb{Z}\OpG}(\underline{\mathbb{Z}},{\rm Pic}) \to \H^2(\OpG,\kt)=\Ext^2_{\mathbb{Z}\OpG}(\underline{\mathbb{Z}},\kt).
$$

\section{Topos cohomology and endotrivial modules}

It is possible to use topos cohomology to demonstrate that Balmer's and Grodal's formulas are naturally identified. It also allows a new characterization of $T_k(G,P)$ by the Picard group of category representations. This section is also based on \cite[Expos\'e V]{SGA4} and \cite[Chapter 21]{St}.

\subsection{Topos cohomology} Topoi are reagarded as generalized spaces. Since a topos can be realized as sheaf categories over different sites, it makes more sense to study cohomology of topoi than that of sites. Suppose $\mathscr{E}$ is a topos, and $\mathfrak{e}$ is a terminal object (which always exists). Let $Ab(\mathscr{E})$ be the category of abelian objects in $\mathscr{E}$, which is abelian and has enough injective objects. Then there is a functor
$$
\Gamma:=\Nat(\mathfrak{e},-): Ab(\mathscr{E}) \to \Ab,
$$
whose derived functors are written as $\H^i(\mathscr{E},-), \forall i\ge 0$. For any object $\frakM\in Ab(\mathscr{E})$, $\H^i(\mathscr{E},\frakM)$ is called the $i$th cohomology of $\mathscr{E}$ with coefficients in $\frakM$. If $\mathscr{E}=\Sh(\C)$ for some site $\C$, then $Ab(\Sh(\C))=\Sh(\C,\Ab)$ and the cohomology group is also expressed as $\H^i(\C,\frakM)$.

We first note that, considering the presheaf topos $\PSh(\calC)$ and an abelian object $\frakM$ in it (that is, a presheaf of abelian groups on $\calC$), we have
$$
\H^i(\PSh(\calC), \frakM)=\H^i(\calC, \frakM).
$$
Thus Grodal's formula $\H^1(\calO_p^{\circ}(G), \kt)$ can be interpreted within topos theory. 

Meanwhile in topos cohomology there is a canonical isomorphism with \v{C}ech cohomology $\H^1(\Sh(\OG), \frakM)=\H^1(\OG,\frakM)\cong\cH^1(G/G,\frakM)$, for any abelian sheaf $\frakM$. The later is identified with $\cH^1(\calU_{G/G},\frakM)$, by Corollary 4.1.4.

\subsection{Leray Spectral sequence} Recall that the inclusion induces a continuous and cocontinuous functor between sites
$$
\iota : (\calO^{\circ}_p(G),\calJ_{\min}) \to (\calO(G),\calJ_{\rm sipp}),
$$
which in turn leads to a geometric morphism between topoi (Corollary 4.3.1)
$$
\mathbb{I}: \PSh(\calO^{\circ}_p(G)) \to \Sh(\OG).
$$
For each $\frakM \in \PSh(\calO^{\circ}_p(G),\Ab)$, there is a Leray spectral sequence \cite[Lemma 21.14.5]{St}
$$
\H^i(\OG, ({\rm R}^j\mathbb{I}_*)(\frakM)) \Rightarrow \H^{i+j}(\calO^{\circ}_p(G), \frakM),
$$
where $\mathbb{I}_*=RK_{\iota}$ and ${\rm R}^j\mathbb{I}_*$ is the $j$-th derived functor. In this spectral sequence, the horizontal edge map reads as follows
$$
\H^i(\OG, \mathbb{I}_*(\frakM)){\buildrel{e}\over{\to}}\H^i(\calO^{\circ}_p(G), \frakM).
$$

To better understand the horizontal edge map, we invoke the Comparison Lemma. We already knew that $(\calO_p(G),\calJ_{\min})$ is dense in $\OG=(\calO(G),\calJ_{\rm sipp})$ and that $\Sh(\OG) \simeq \PSh(\calO_p(G))$. Thus we may equally turn to the Leray spectral sequence for the following functor, still denoted by $\iota$,
$$
\iota : (\calO^{\circ}_p(G),\calJ_{\min}) \to (\calO_p(G),\calJ_{\min}).
$$
and the induced geometric morphism between the (presheaf) topoi
$$
\mathbb{I}: \PSh(\OpG) \to \PSh(\calO_p(G)).
$$
For each $\frakM \in \PSh(\OpG,\Ab)$, there is a Leray spectral sequence
$$
\H^i(\calO_p(G), ({\rm R}^j\mathbb{I}_*)(\frakM)) \Rightarrow \H^{i+j}(\calO^{\circ}_p(G), \frakM),
$$
where $\mathbb{I}_*=RK_{\iota}$ and ${\rm R}^j\mathbb{I}_*$ is the $j$-th derived functor. In this spectral sequence, the horizontal edge map is
$$
\H^i(\calO_p(G), \mathbb{I}_*(\frakM)){\buildrel{e}\over{\to}}\H^i(\calO^{\circ}_p(G), \frakM),
$$
which is the composition of
$$
\H^i(\calO_p(G)), \mathbb{I}_*(\frakM))\to\H^i(\calO^{\circ}_p(G), \Res_{\iota}{\mathbb{I}}_*(\frakM))\to\H^i(\calO^{\circ}_p(G), \frakM),
$$
entirely analogous to the case of the Lydon-Hochschild-Serre spectral sequence for a group extension. The map on the left is induced by the restriction $Res_{\iota}$. While the other one is given by the natural map $Res_{\iota}{\mathbb{I}}_*(\frakN)\cong\frakN \to \frakN$ from adjunction. In our case, since $Res_{\iota}{\mathbb{I}}_*(\frakN)$ is the restriction of the right Kan extension of $\frakN$ along the fully faithful $\iota : \OpG\to\calO_p(G)$, we have $Res_{\iota}{\mathbb{I}}_*(\frakN)\cong\frakN$. It implies that the horizontal edge map
$$
\H^i(\calO_p(G)), \mathbb{I}_*(\frakM)){\buildrel{e}\over{\to}}\H^i(\calO^{\circ}_p(G), \frakM)
$$
is exactly the one induced by $Res_{\iota}$.
 
\begin{proposition}
The preceding horizontal edge map is an isomorphism. 
\end{proposition}

\begin{proof} It may be expressed as 
$$
\Ext^i_{\mathbb{Z}\calO_p(G)}(\underline{\mathbb{Z}},\mathbb{I}_*(\frakM)){\buildrel{e}\over{\to}}\Ext^i_{\mathbb{Z}\OpG}(\underline{\mathbb{Z}},\frakM),
$$
as the map induced by $Res_{\iota}: \PSh(\calO_p(G), \Ab) \to \PSh(\OpG, \Ab)$. Proceed as in the proof of Corollary 5.1.3, we are done.
\end{proof}

\begin{remark}
The horizontal edge map is some sort of ``inflation'' if we compare it with the spectral sequence for group extension $\pi: G \to G/N$, where $N$ is a normal subgroup of $G$. When we regard these two groups as categories, the functor $\pi$ induces a topos morphism $\PSh(G)\to\PSh(G/N)$ \cite[IV 4.5]{SGA4}. The topos cohomology of $\PSh(G)$ is just the usual group cohomology. The Leray spectral sequence for this topos morphism is the Lydon-Hochschild-Serre spectral sequence for the group extension $1\to N \to G \to G/N \to 1$. 
\end{remark}

\subsection{The five term exact sequence} In its original form, we know the horizontal edge map
$$
\H^i(\OG, \mathbb{I}_*(\frakM)){\buildrel{e}\over{\to}}\H^i(\calO^{\circ}_p(G), \frakM)
$$
is an isomorphism for all $i\ge 0$ and all $\frakM\in\PSh(\OpG, \Ab)$. If we turn to the resulting five term exact sequence 
$$
0\to\H^1(\OG, \mathbb{I}_*(\frakM)){\buildrel{e}\over{\to}}\H^1(\calO^{\circ}_p(G), \frakM){\buildrel{e'}\over{\to}}
$$
$$
\H^1(\OG, ({\rm R}^1\mathbb{I}_*)(\frakM)){\buildrel{d}\over{\to}}\H^2(\OG, \mathbb{I}_*(\frakM)){\buildrel{e}\over{\to}}\H^2(\calO^{\circ}_p(G), \frakM),
$$
then we obtain the following observation.

\begin{corollary} For each $\frakM\in\PSh(\OpG, \Ab)$, 
$$\H^1(\OG, ({\rm R}^1\mathbb{I}_*)(\frakM))=0.$$
\end{corollary}

Particularly the map $e:\H^1(\OG, \mathbb{G}_m)\to\H^1(\calO^{\circ}_p(G), \underline{k^{\times}})$ offers another way to identify Balmer's and Grodal's formulas on $T_k(G,P)$.

\subsection{$T_k(G,P)$ as a Picard group} By Corollaries 4.1.4 and 5.1.3, we have seen that
$$
T_k(G,P)\cong\cH^1(\calU_{G/G},\GM)\cong\cH^1(G/G,\GM).
$$
Let $\frakO$ be the presheaf of rings such that $\frakO(G/H)=k$ if $p\bigm{|}|H|$, and zero otherwise. Similar to the proof of Proposition 3.2.3, $\frakO=\I_*(\underline{k})$ turns out to be a sheaf. Let $\frakO^{\times}$ be the subsheaf of abelian groups defined by $\frakO^{\times}(G/H):=\frakO(G/H)^{\times}=k^{\times}$, the group of invertible elements of $\frakO(G/H)$, $\forall G/H$. Then $\frakO^{\times}=\GM$. By \cite[Lemma 21.6.1]{St} there is an isomorphism of abelian groups
$$
\cH^1(G/G,\GM)\cong\H^1(\OG,\GM)\cong\Pic(\frakO),
$$
where $\Pic(\frakO)$ is the Picard group \cite[Section 18.32]{St} of the ringed site $(\OG,\frakO)$, consisting of isomorphism classes of invertible $\frakO$-modules (the tensor product between invertible modules gives rise to the group structure). Meanwhile in \cite[Theorem 4.2.1]{WX}, we had an equivalence of ringed sites (or ringed topoi)
$$
(\OG,\frakO)\simeq(\calO_p(G),\frakO|_{\calO_p(G)}).
$$
By the definition of $\frakO$, we continue to deduce another equivalence of ringed sites (both carrying the minimal topology)
$$
(\calO_p(G),\frakO|_{\calO_p(G)})\simeq(\calO^{\circ}_p(G),\underline{k}).
$$
Therefore
$$
\Pic(\frakO)\cong\Pic(\underline{k}).
$$

In fact, we can arrive at the same observation from Grodal's formula to get $T_k(G,P)\cong\H^1(\calO^{\circ}_p(G),\underline{k^{\times}})\cong\Pic(\underline{k})$, regarding $\H^1(\calO^{\circ}_p(G),\underline{k^{\times}})$ as a presheaf topos cohomology group, replacing $\H^1(\OG,\GM)$ above.

\begin{theorem} There is an isomorphism of abelian groups $T_k(G,P)\cong \Pic(\underline{k})$, the Picard group of the ringed site $(\calO^{\circ}_p(G), \underline{k})$.

The group $\Pic(\underline{k})$ consists of the isomorphism classes of presheaves $\frakF$ satisfying the conditions that $\dim_k\frakF(G/Q)=1$ and that  $\frakF(G/Q\to G/Q')$ is an isomorphism of $k$-vector spaces, for every morphism $G/Q \to G/Q'$ of $\calO^{\circ}_p(G)$.
\end{theorem}

\begin{proof} We only need to prove the second statement. It is relatively easier to understand invertible modules in module categories than in stable module categories. The elements of $\Pic(\k)$ are represented by the isomorphism classes of invertible presheaves on $\calO^{\circ}_p(G)$, with values in finite-dimensional $k$-vector spaces (since invertible modules are of finite presentation). The group operation is given by the pointwise tensor product between two presheaves. 

If $\frakF$ is an invertible presheaf on $\calO^{\circ}_p(G)$, then $\frakF(G/Q)$ is an invertible $kN_G(Q)/Q$-module for every non-trivial $p$-subgroup $Q$ of $G$, and $\frakF(G/Q\to G/Q')\ne 0$ on any morphism.  It implies that $\dim_k\frakF(G/Q)=1, \forall Q$, and that $\frakF(G/Q\to G/Q')$ is an isomorphism on every morphism $G/Q\to G/Q'$ of $\calO^{\circ}_p(G)$. Given a presheaf $\frakF$ on $\calO^{\circ}_p(G)$, satisfying the two conditions, we are able to define a presheaf $\frakF^*$ on $\calO^{\circ}_p(G)$ such that
$$
\frakF^*(G/Q) := \frakF(G/Q)^*
$$
and such that (as a $k$-linear map) $\frakF^*(G/Q \to G/Q')$ is the $k$-dual of the inverse map $\frakF(G/Q\to G/Q')^{-1}$. More explicitly, if $a_{Q'}$ is any element of $\frakF(G/Q')$ and $\frakF(G/Q\to G/Q')(a_{Q'})=a_Q$, then
$$
\frakF^*(G/Q \to G/Q')(a^*_{Q'}) := a^*_Q.
$$
Here $a^*_{Q'}$ (or $a_Q^*$) is the dual element of $a_{Q'}$ (or $a_Q$). One can readily verify that $\frakF\otimes\frakF^*\cong\underline{k}$. 
\end{proof}

\begin{remark} 
\begin{enumerate}
\item By Remark 3.2.8, every $kG$-module $M$ gives rise to a fixed-point presheaf on $\calO^{\circ}_p(G)$ by $\frakF_M(G/Q)=M^Q$, for all non-trivial $p$-subgroups $Q$. If $M$ is 1-dimensional (hence endotrivial), then $M^Q=M$. From $\frakF_{M\otimes N}\cong\frakF_M\otimes\frakF_N$, under the circumstance $\frakF_M$ becomes a (constant-valued but not necessarily constant) invertible presheaf on $\calO^{\circ}_p(G)$. This establishes an injection from the set of the isomorphism classes of 1-dimensional $kG$-modules to that of the Sylow-trivial modules, that is, 
$$
\Hom(G,k^{\times}) \hookrightarrow \Pic(\underline{k})\cong T_k(G,P).
$$

\item The group $\Pic(\frakO)\cong\Pic(\underline{k})$ consists of the isomorphism classes of presheaves $\frakG$ on $\calO_p(G)$ such that $\dim_k\frakG(G/Q)=1, \forall Q\ne 1$ and $\frakG(G/1)=0$, and such that  $\frakG(G/Q\to G/Q')$ is an isomorphism of $k$-vector spaces, for every morphism $G/Q \to G/Q'$ of $\calO_p(G)$, whenever $Q\ne 1$.

\item One may construct an explicit isomorphism between the Picard group and the degree one topos cohomology group, see \cite[Lemma 21.6.1]{St}.
\end{enumerate}
\end{remark}

\end{document}